\documentclass{article}
\usepackage{amsmath}
\usepackage{amsthm}
\usepackage{color}
\usepackage{appendix}
\usepackage{booktabs}
\usepackage[justification=centerfirst]{caption}
\usepackage{hyperref}
\usepackage{hypcap}

\begin{document}

\title{Compounding Doubly Affine Matrices}
\author{Adam Rogers, Ian Cameron and Peter Loly \\ Department of Physics and Astronomy, University of Manitoba,\\Winnipeg, Manitoba, Canada R3T 2N2}
\maketitle

\begin{abstract}

Weighted sums of left and right hand Kronecker products of Integer Sequence
Doubly Affine (ISDA) as well as Generalized Arithmetic Progression Doubly Affine (GAPDA) arrays are used to generate larger ISDA arrays of multiplicative order (\textit{compound squares}) from pairs of smaller ones.

In two dimensions we find general expressions for the eigenvalues (EVs) and
singular values (SVs) of the larger arrays in terms of the EVs and SVs of
their constituent matrices, leading to a simple result for the rank of
these highly singular compound matrices. Since the critical property of the
smaller constituent matrices involves only identical row and column sums
(often called semi-magic), the eigenvalue and singular value results can be
applied to both magic squares and Latin squares. Additionally, the
compounding process works in arbitrary dimensions due to the generality of the
Kronecker product, providing a simple method to generate large order ISDA
cubes and hypercubes.

The first examples of compound magic squares are found in manuscripts that
date back to the 10th century CE, and other representative applications are
outlined through judicious examples.


\end{abstract}

\section{Introduction}
\label{sec1}

Magic squares are $n \times n$ arrays ($n \geq 3$), comprised of the distinct elements $0..(n^2-1)$. The elements appear only once in a given square, called full (or exact) cover, and are arranged such that every row, column, and both main diagonals (RCD)
add to the same value, often called the magic sum $S_{n,k}$. Semi-magic squares are defined in analogy to magic squares but lack one or both of the diagonal sums. When the restriction to unique elements is lifted, the set of semi-magic squares includes Latin squares, which use the complete sequence $0..(n-1)$ with ($n \geq 2$) in each row of the matrix, arranged such that the rows and columns (but not necessarily the diagonals) sum to the same constant value. We categorize the magic, semi-magic and Latin squares as examples of ``integer sequence doubly affine'' (ISDA) matrices \cite{powers}. An integer parameter $k$ serves to distinguish Latin arrays ($k=1,$ $n\geq
2$), and magic squares ($k=2,n\geq3$) \cite{Drury} with elements $0..(n^{k}-1)$, so that the line sum, for order $n$ in both cases has a particularly simple expression:
\begin{equation}
S_{n,k}=\frac{n}{2}\left(  n^{k}-1\right).
\label{magicSum}
\end{equation}

In general, we can extend the definition of an ISDA matrix to an ISDA array in an arbitrary number of dimensions $D$ using larger integer sequences. For example, magic cubes are cubic arrays of integers $0..(n^3-1)$ with sums along each of the rows, columns, pillars and body diagonals that sum to a constant value. Generally, we write the line sum constant for a magic object containing sequential integer elements $0..(n^k-1)$, such that $n \times n$ Latin arrays (squares, cubes, hypercubes and higher dimensional structures) have $k=1$ ($S_{n,1}$), magic squares have $k=2$, magic
cubes have $k=3$, etc.

An algorithm for combining order $n$ and $m$ ISDA ``seed'' matrices to result in a larger order $mn$ ``compound'' product was developed by Loly and Chan \cite{chan}. That work focused on magic squares in particular. In later work, Loly et al \cite{LCTS,LCTRS} studied the eigenspectra $\lambda_i$ (which may be real, complex or simply imaginary) as well as the singular value $\sigma_i^2$ spectra of magic squares. Recently Cameron, Rogers \& Loly \cite{ian2},\cite{CRL} developed the
concept of singular value clans to characterize the singular value sets which
identify magic or Latin squares related by various permutations (see also
Staab et al \cite{Staab}). Our work included new measures flowing from the
SVs ($\sigma_{i}$), namely Shannon entropy measures and especially an integer index:
\begin{equation}
R=\sum_{i=2}^{n} \sigma_{i}^{4}
\end{equation}
which sums the fourth powers of the SVs after the first linesum SV, and
identifies ``clans'' of related squares. For example, $144=3\times48$ of the
$880$ distinct order four magic squares fall into just three clans
\cite{LCTRS}\cite{LCTS}, each shared by $16$ members of Dudeney's
\cite{Dudeney} Group I (pandiagonal), II (semi-bent, semi-pandiagonal) and III
(regular, semi-pandiagonal). Within each clan and Dudeney Group the $16$ members are
all related by row and column permutations (Conway \cite{Conway}).

Another issue concerns the propagation of common special properties of the
seed squares under compounding. Chan and Loly \cite{chan} gave an argument
for the propagation of pandiagonality under aggregated compounding, and
Eggermont \cite{Egger} proved the propagation of those and a larger set of
properties, all for the magic square aggregated case. Eggermont's thesis
\cite{Egger} is a valuable resource for the aggregated
compounding of magic squares. Since Latin squares may not be strongly
diagonal, we show explicitly that compounding a strongly diagonal Latin square
with another which is not diagonal results in compound Latin squares with the lesser property.

In this work we describe an analytical process using Kronecker products \cite{Basdevant} to systemmatize the construction of compound ISDA arrays. In Section \ref{sec:historical} we show examples of historical squares which can be described in terms of the compounding operation. In Section \ref{sec:compoundingISDA} we describe the construction method for building compound arrays in terms of smaller ISDA arrays, and then generalize the approach to use more general arrays that include gaps in the integer sequence in Section \ref{sec:compoundingGAPDA}. In two dimensions we derive the eigenproperties and the singular value decomposition (SVD) of the compound matrices in terms of the lower order constituent matrices. Our study of the SVD also determines the rank of the compound matrices. We emphasize that this general process can be used to construct any order ISDA arrays such as Latin ($k=1$) or magic ($k \geq 2$) arrays in arbitrary dimension.

\section{Historical Compound Squares}
\label{sec:historical}

Compound magic squares have been rediscovered many times. A particularly
simple method of construction produces a larger magic square of
multiplicative order from two smaller ``seed'' squares that are appropriately
tiled and incremented \cite{Andrews, Schubert}. To illustrate the compounding procedure, consider the smallest magic square ($n=3$), which is unique if we neglect the 8 geometric symmetries of rotation and reflection,
\begin{equation}
\mathbf{M}_{3}=\left[
\begin{array}
[c]{ccc}%
3 & 8 & 1\\
2 & 4 & 6\\
7 & 0 & 5
\end{array}
\right].
\label{M3}
\end{equation}
This ISDA matrix is called the Lo Shu square. Note that in this work we are using the sequence $0..(n^2-1)$ for mathematical convenience, however magic squares are historically written over $1..n^2$. This rank $r=3$ square has eigenspectra $\lambda_i=\{12,+2i\sqrt{3},-2i\sqrt{3}\}$ and singular values $\sigma_i^2=\{144,48,12\}$ with index $R=2304$. The singular values are simply $\sigma_i$, but in this work we will give the squared values $\sigma_i^2$ unless otherwise indicated.

The geometric arrangement of the ordered integer sequence $\{0,1,...,8\}$ has a criss-cross pattern in the Lo Shu that is reflected in the structure of the historical compound squares, both in the placement of the subsquares and within the subsquares themselves. The first type of compound square was originally written prior to $1000$ AD \cite{SesianoCM}:
\begin{equation}
\mathbf{C}_{\text{A}}=\left[
\begin{array}
[c]{ccc}%
\begin{array}
[c]{ccc}%
30 & 35 & 28\\
29 & 31 & 33\\
34 & 27 & 32
\end{array}
&
\begin{array}
[c]{ccc}%
75 & 80 & 73\\
74 & 76 & 78\\
79 & 72 & 77
\end{array}
&
\begin{array}
[c]{ccc}%
12 & 17 & 10\\
11 & 13 & 15\\
16 & 9 & 14
\end{array}
\\%
\begin{array}
[c]{c}%
\;\\%
\begin{array}
[c]{ccc}%
21 & 26 & 19\\
20 & 22 & 24\\
25 & 18 & 23
\end{array}
\end{array}
&
\begin{array}
[c]{c}%
\;\\%
\begin{array}
[c]{ccc}%
39 & 44 & 37\\
38 & 40 & 42\\
43 & 36 & 41
\end{array}
\end{array}
&
\begin{array}
[c]{c}%
\;\\%
\begin{array}
[c]{ccc}%
57 & 62 & 55\\
56 & 58 & 60\\
61 & 54 & 59
\end{array}
\end{array}
\\%
\begin{array}
[c]{c}%
\;\\%
\begin{array}
[c]{ccc}%
66 & 71 & 64\\
65 & 67 & 69\\
70 & 63 & 68
\end{array}
\end{array}
&
\begin{array}
[c]{c}%
\;\\%
\begin{array}
[c]{ccc}%
\;\textcolor{red}{3} & \;\textcolor{red}{8} & \;\textcolor{red}{1}\\
\;\textcolor{red}{2} & \;\textcolor{red}{4} & \;\textcolor{red}{6}\\
\;\textcolor{red}{7} & \;\textcolor{red}{0} & \;\textcolor{red}{5}
\end{array}
\end{array}
&
\begin{array}
[c]{c}%
\;\\%
\begin{array}
[c]{ccc}%
48 & 53 & 46\\
47 & 49 & 51\\
52 & 45 & 50
\end{array}
\end{array}
\end{array}
\right].  \label{Arabic}%
\end{equation}
We label this square $\mathbf{C}$ to show it is a compound square which we call ``aggregated'', and label with the subscript $\text{A}$, as it consists of augmentations of the $3\times3$ Lo Shu subsquares it contains. For example, the bottom center block (shown in red text) is just the Lo Shu, with the others incremented by factors of $9$ and added to the corresponding positions of the integers in the compound square \cite{SesianoAbuW, ThompsonWH}. This basic aggregated compounding scheme can be used with pairs of ISDA arrays that lead to larger arrays of "multiplicative order", e.g. $3 \times 4=12$, $3 \times 5=15$, $4 \times 4=16$, etc. The subsquare structure is emphasized by the blank medial row and column, a feature which will be used throughout this work.

Simply using the Lo Shu, it is possible to form a second type of compound square, again in standard form, which contains magic $3 \times 3$ subsquares with non-sequential elements:
\begin{equation}
\mathbf{C}_{\text{D}}=\left[
\begin{array}
[c]{ccc}%
\begin{array}
[c]{ccc}%
30 & 75 & 12\\
21 & 39 & 57\\
66 & \textcolor{red}{3} & 48
\end{array}
&
\begin{array}
[c]{ccc}%
35 & 80 & 17\\
26 & 44 & 62\\
71 & \textcolor{red}{8} & 53
\end{array}
&
\begin{array}
[c]{ccc}%
28 & 73 & 10\\
19 & 37 & 55\\
64 & \textcolor{red}{1} & 46
\end{array}
\\%
\begin{array}
[c]{c}%
\;\\%
\begin{array}
[c]{ccc}%
29 & 74 & 11\\
20 & 38 & 56\\
65 & \textcolor{red}{2} & 47
\end{array}
\end{array}
&
\begin{array}
[c]{c}%
\;\\%
\begin{array}
[c]{ccc}%
31 & 76 & 13\\
22 & 40 & 58\\
67 & \textcolor{red}{4} & 49
\end{array}
\end{array}
&
\begin{array}
[c]{c}%
\;\\%
\begin{array}
[c]{ccc}%
33 & 78 & 15\\
24 & 42 & 60\\
69 & \textcolor{red}{6} & 51
\end{array}
\end{array}
\\%
\begin{array}
[c]{c}%
\;\\%
\begin{array}
[c]{ccc}%
34 & 79 & 16\\
25 & 43 & 61\\
70 & \textcolor{red}{7} & 52
\end{array}
\end{array}
&
\begin{array}
[c]{c}%
\;\\%
\begin{array}
[c]{ccc}%
27 & 72 & 9\\
18 & 36 & 54\\
63 & \textcolor{red}{0} & 45
\end{array}
\end{array}
&
\begin{array}
[c]{c}%
\;\\%
\begin{array}
[c]{ccc}%
32 & 77 & 14\\
23 & 41 & 59\\
68 & \textcolor{red}{5} & 50
\end{array}
\end{array}
\end{array}
\right]  . \label{YH}%
\end{equation}
Here the subscript D in $\mathbf{C}_{\text{D }}$ refers to the dispersed positions of sequential integers in different subsquares. This second type of compound square dates to 1275 AD \cite{Lam, LiYen, Cam60, Cammann, Swetz}. Unlike the aggregated compound square, there has been apparently no development of an underlying construction method to produce large order dispersed squares from single ISDA arrays or from pairs of ISDA arrays.  It is worth noting that the two squares are perfect two dimensional multiway shuffles of each other (a particular interleaving permutation that reorders the row and column sequence ($012\cdot 345\cdot 678$) to ($036\cdot 147\cdot 258$) by splitting the original sequence into 3 groups and selecting a row or column number in corresponding position from each group) and so there is no fundamental change to the structure of the squares, ie the individual elements in the rows and columns maintain the same relationship to one another. Indeed, the elements that occupied the bottom center subsquare locations of the previous square $\mathbf{C}_{\text{A}}$ appear now all together condensed into the bottom central subsquare of $\mathbf{C}_{\text{D}}$ with no change to their earlier Lo Shu pattern.  (Note also that the diagonal elements have simply been rearranged). Both of these rank $r=5$ compound squares have the same eigenspectra $\lambda_i=\{360,+54i\sqrt{6},-54i\sqrt{6},+6i\sqrt{6},-6i\sqrt{6},0,0,0,0\}$ and singular values $\sigma_i^2=\{360^2,\allowbreak 3\times108^2,3\times54^2,3\times12^2,3\times6^2,0,0,0,0\}$ with index $R=1,301,165,856$.

A more recent construction is the Korean Koo-Soo-Ryak square, which can be traced to Choi Seok-Jeong (1645-1715) \cite{KSR}:

\begin{equation}
\mathbf{C}_{\text{C}}=\left[
\begin{array}
[c]{ccc}%
\begin{array}
[c]{ccc}%
 36 & 47 & 28 \\
 29 & 37 & 45 \\
 46 & 27 & 38
\end{array}
&
\begin{array}
[c]{ccc}%
69 & 80 & 61 \\
62 & 70 & 78 \\
79 & 60 & 71
\end{array}
&
\begin{array}
[c]{ccc}%
12 & 23 & \textcolor{red}{4}\\
  \textcolor{red}{5} & 13 & 21\\
22 &   \textcolor{red}{3} & 14
\end{array}
\\%
\begin{array}
[c]{c}%
\;\\%
\begin{array}
[c]{ccc}%
15 & 26 & \textcolor{red}{7} \\
  \textcolor{red}{8} & 16 & 24\\
25 &   \textcolor{red}{6} & 17
\end{array}
\end{array}
&
\begin{array}
[c]{c}%
\;\\%
\begin{array}
[c]{ccc}%
 39 & 50 & 31\\
 32 & 40 & 48\\
 49 & 30 & 41
\end{array}
\end{array}
&
\begin{array}
[c]{c}%
\;\\%
\begin{array}
[c]{ccc}%
63 & 74 & 55\\
56 & 64 & 72\\
73 & 54 & 65
\end{array}
\end{array}
\\%
\begin{array}
[c]{c}%
\;\\%
\begin{array}
[c]{ccc}%
66 & 77 & 58\\
59 & 67 & 75\\
76 & 57 & 68
\end{array}
\end{array}
&
\begin{array}
[c]{c}%
\;\\%
\begin{array}
[c]{ccc}%
  9 & 20 & \textcolor{red}{1}\\
  \textcolor{red}{2} & 10 & 18\\
19 &   \textcolor{red}{0} & 11
\end{array}
\end{array}
&
\begin{array}
[c]{c}%
\;\\%
\begin{array}
[c]{ccc}%
42 & 53 & 34\\
35 & 43 & 51\\
52 & 33 & 44
\end{array}
\end{array}
\end{array}
\right]  . \label{KSR}
\end{equation}
The Koo-Soo-Ryak square, though quite different in its arrangement of actual numbers, is just a more complicated example of compounding, but cannot be reproduced using a simple tiling method from ISDA squares. Therefore we will generalize our compounding formulae to use generalized arithmetic progression doubly affine arrays (GAPDAs). These are arrays which contain gaps in the integer sequence of their elements, but when compounded produce larger order ISDA squares with full cover of the integers $0..(n^2-1)$. In this instance the two GAPDA's are derived from the generalized arithmetic progressions (GAPs): $\{0,1,2,9,10,11,18,19,20\}$ and $\{0,3,6,27,30,33,54,57,60\}$ by arranging each of them in the Lo Shu configuration. Note the appearance of all the elements of the first array in the zero position subsquare of $\mathbf{C}_\text{C}$ and the elements of the second array in the zero positions of all the subsquares of $\mathbf{C}_\text{C}$.

Finally, we note that the three order $9$ squares $\mathbf{C}_\text{A}$, $\mathbf{C}_\text{D}$ and $\mathbf{C}_\text{C}$ were also given in the work of Frierson \cite{Frierson}. In that work, an algebraic parameterization of the Lo Shu led to a set of six order $9$ compound squares. Here, we will derive these six squares using our analytical compounding formulae and the appropriate seed squares.

\section{Compounding ISDA arrays}
\label{sec:compoundingISDA}

The compounding method using ISDA seeds relies on the appropriate tiling and incrementing of the constituent arrays. In order to describe this process mathematically, we make use of particular weighted sums of Kronecker products. For simplicity we will develop the method in two dimensions, entirely using matrices. However, we stress that our formulae are valid for arrays of arbitrary dimensions.

Consider a matrix $\mathbf{E}_{m}$ of size $m \times m$ with each element equal
to $1$, and an arbitrary ISDA matrix $\mathbf{B}_{n}$ of size $n \times n$, containing all integers $0..(n^2-1)$. Define the Kronecker product \cite{Basdevant} between these two matrices as follows:
\begin{equation}
\mathbf{\beta}_{n}=\mathbf{E}_{m}\otimes\mathbf{B}_{n}=
\left[
\begin{array}
[c]{ccc}
\mathbf{B}_{n} & ... & \mathbf{B}_{n}\\
\vdots & \ddots & \vdots\\
\mathbf{B}_{n} & ... & \mathbf{B}_{n}
\end{array}
\right].
\label{beta}
\end{equation}
This operation produces an $mn\times mn$ matrix $\mathbf{\beta}_{n}$ which consists of $m$ tiled copies of $\mathbf{B}_{n}$. Let us now consider a second Kronecker product, between an $m \times m$ matrix $\mathbf{A}_{m}$ and an $n \times n$ matrix $\mathbf{E}_{n}$, consisting of $1$ in every entry:
\begin{equation}
\text{ }
\mathbf{\alpha}_{m}=\mathbf{A}_{m}\otimes\mathbf{E}_{n}=
\left[
\begin{array}
[c]{ccc}
a_{11}\mathbf{E}_{n} & ... & a_{1m}\mathbf{E}_{n} \\
\vdots & \ddots & \vdots \\
a_{m1}\mathbf{E}_{n} & ... & a_{mm}\mathbf{E}_{n}
\end{array}
\right].
\label{alpha}
\end{equation}
The resulting $mn \times mn$ matrix $\mathbf{\alpha}_{m}$ is composed of square $n\times n$ tiles of each element of $\mathbf{A}_{m}$. With these simple definitions, we are now equipped to explore the compounding procedure.

To illustrate a simple case of the compounding procedure for ISDA square matrices of different orders. We begin with the smallest Latin square of two symbols, chosen here as $0,1$:
\begin{equation}
\mathbf{L}_2=\left[
\begin{array}
[c]{cc}
0 & 1\\
1 & 0
\end{array}
\right],
\end{equation}
where the rows and columns share a common row-column (RC) linesum. Let us compound this square with the third order Latin square,
\begin{equation}
\mathbf{L}_3=\left[
\begin{array}
[c]{ccc}
0 & 1 & 2\\
1 & 2 & 0\\
2 & 0 & 1
\end{array}
\right].
\end{equation}
with $n=3$. This rank $r=3$ square has eigenvalues $\mu_1=3$, $\mu_2=\sqrt{3}$ and $\mu_3=-\mu_2$ and singular values $\sigma_i^2=\{9,3,3\}$ with index $R=99\textit{(=81+9+9)}$. We will take $m=2$, such that $\mathbf{A}_m=\mathbf{L}_2$ and $\mathbf{B}_n=\mathbf{L}_3$, with ($\mathbf{A}_m = \mathbf{L}_2$, $\mathbf{B}_n = \mathbf{L}_3$, $m=2$, $n=3$ and $k=1$). In this case, we find
\begin{equation}
\mathbf{\alpha}_{2} = \mathbf{L}_2 \otimes \mathbf{E}_{3} = \left[
\begin{array}
[c]{ccccccc}
0 & 0 & 0 &  & 1 & 1 & 1 \\
0 & 0 & 0 &  & 1 & 1 & 1 \\
0 & 0 & 0 &  & 1 & 1 & 1 \\
  &   &   &  &   &   &   \\
1 & 1 & 1 &  & 0 & 0 & 0 \\
1 & 1 & 1 &  & 0 & 0 & 0 \\
1 & 1 & 1 &  & 0 & 0 & 0
\end{array}
\right],
\end{equation}
\begin{equation}
\mathbf{\beta}_{3} = \mathbf{E}_{2} \otimes \mathbf{L}_3 = \left[
\begin{array}
[c]{ccccccc}
0 & 1 & 2 &  & 0 & 1 & 2 \\
1 & 2 & 0 &  & 1 & 2 & 0 \\
2 & 0 & 1 &  & 2 & 0 & 1 \\
  &   &   &  &   &   &   \\
0 & 1 & 2 &  & 0 & 1 & 2 \\
1 & 2 & 0 &  & 1 & 2 & 0 \\
2 & 0 & 1 &  & 2 & 0 & 1
\end{array}
\right].
\end{equation}
We can then add these matrices in combination to produce the aggregated compound Latin square by multiplying $\mathbf{\alpha}_2$ by $n^k=3$ and adding it to $\mathbf{\beta}_3$ to give
\begin{equation}
\label{isda0}
\mathbf{L}_{6A}=\left[
\begin{array}
[c]{ccccccc}%
\textcolor{blue}0 & 1 & 2 & & \textcolor{blue}3 & 4 & 5\\
1 & 2 & 0 & & 4 & 5 & 3\\
2 & 0 & 1 & & 5 & 3 & 4\\
  &   &   & &   &   &  \\
\textcolor{blue}3 & 4 & 5 & & \textcolor{red}0 & \textcolor{red}1 & \textcolor{red}2\\
4 & 5 & 3 & & \textcolor{red}1 & \textcolor{red}2 & \textcolor{red}0\\
5 & 3 & 4 & & \textcolor{red}2 & \textcolor{red}0 & \textcolor{red}1
\end{array}
\right]
\end{equation}
where we have emphasized its structure with colored blocks. The aggregated square has non-vanishing eigenvalues $\lambda_1=15$, $\lambda_2=-9$, $\lambda_3=2\sqrt{3}$ and $\lambda_4=-\lambda_3$, $r=4$, $\sigma_i^2=\{225,81,12,12,0,0\}$, and index $R=6849$.

While we have used a simple Latin square ($k=1$) to clearly illustrate the method, the basic tiling process works in general. Given the constituent ISDA arrays $\mathbf{A}_{m}$ and $\mathbf{B}_{n}$, we define the general aggregated compound matrix
\begin{equation}
\mathbf{C}_{A}=n^{k}\mathbf{\alpha}_{m}+\mathbf{\beta}_{n}.
\label{cmpA}
\end{equation}
with line sum
\begin{equation}
S_{C,k}=n^{k+1}S_{m,k}+mS_{n,k}=\frac{mn}{2}\left[\left(mn\right)^{k}-1\right]
\label{cmpLS}
\end{equation}
Note that when $\mathbf{A}_{m}=\mathbf{B}_{n}=\mathbf{M}_3$ where $\mathbf{M}_3$ is the Lo Shu magic square (\ref{M3}), we find that (\ref{cmpA}) produces the aggregated compound square in (\ref{Arabic}) with $m=n=3$ and $k=2$.

It is possible to construct the dispersed type of compound square (\ref{YH}) from $\mathbf{A}_m$ and $\mathbf{B}_n$ using an appropriate combination of $\mathbf{\alpha}$ and $\mathbf{\beta}$. Dispersed squares are generated by
\begin{equation}
\mathbf{C}_{D}=\mathbf{\alpha}_{m}+m^{k}\mathbf{\beta}_{n}.
\label{cmpD}
\end{equation}
The example above produces the dispersed square
\begin{equation}
\label{isda1}
\mathbf{L}_{6D}=\left[
\begin{array}
[c]{ccccccc}%
\textcolor{blue}0 & 2 & 4 & & \textcolor{blue}1 & 3 & 5\\
2 & 4 & 0 & & 3 & 5 & 1\\
4 & 0 & 2 & & 5 & 1 & 3\\
  &   &   & &   &   &  \\
\textcolor{blue}1 & 3 & 5 & & \textcolor{red}0 & \textcolor{red}2 & \textcolor{red}4\\
3 & 5 & 1 & & \textcolor{red}2 & \textcolor{red}4 & \textcolor{red}0\\
5 & 1 & 3 & & \textcolor{red}4 & \textcolor{red}0 & \textcolor{red}2
\end{array}
\right],
\end{equation}
with the block structure emphasized with text color. The dispersed square has $\lambda_1=15$, $\lambda_2=-3$, $\lambda_3=4\sqrt{3}$ and $\lambda_4=-\lambda_3$, singular values $\sigma_i^2=\{225,48,48,9,0,0\}$, index $R=4689$.

We note that the first (aggregated) type has also been studied by what amounts
to Kronecker products by Derksen, Eggermont \& van den Essen \cite{derksen1}
(see the thesis of Eggermont \cite{Egger} who characterized the pattern as
3-partitioned). The second (dispersed) type compound magic square has not been fully investigated until the present work. These aggregated and distributed types constitute a duet. More recently, we have recognized an additional pair if the seed squares are of different order, for what amounts then to a fundamental quartet. (For Latin squares in general, we are able to find an additional pair for all compound squares of small composite order $mn$, $m \neq n$ and for prime power compounds $p^{q+r}$ where $q+r>2$). The aggregated and dispersed squares of the second duet will now have a different structure by exhibiting a new tiling, and will be perfect shuffles of the first.

We define the ``reverse'' arrays
\begin{equation}
\mathbf{\alpha}_{m}^{\mathcal{R}}=\mathbf{E}_{n} \otimes \mathbf{A}_{m}
\end{equation}
and
\begin{equation}
\mathbf{\beta}_{n}^{\mathcal{R}}=\mathbf{B}_{n} \otimes \mathbf{E}_{m}
\end{equation}
to define the reverse aggregated and dispersed matrices
\begin{equation}
\mathbf{C}_{A}^{\mathcal{R}}=n^{k}\mathbf{\alpha}_{m}^{\mathcal{R}}+\mathbf{\beta}_{n}^{\mathcal{R}}
\label{revCmpA}
\end{equation}
\begin{equation}
\mathbf{C}_{D}^{\mathcal{R}}=\mathbf{\alpha}_{m}^{\mathcal{R}}+m^{k}\mathbf{\beta}_{n}^{\mathcal{R}}.
\label{revCmpD}
\end{equation}
Using the $\mathbf{L}_2$ and $\mathbf{L}_3$ seed squares from above, we find the reverse aggregated and dispersed squares,
\begin{equation}
\label{isda2}
\mathbf{L}_{6A}^{\mathcal{R}}=\left[
\begin{array}
[c]{ccccccc}
0 & 3 & 1 & &  4 & 2 & 5\\
3 & 0 & 4 & &  1 & 5 & 2\\
1 & 4 & 2 & &  5 & 0 & 3\\
   &    &   &  &   &    &   \\
4 & 1 & 5 & &  2 & 3 & 0\\
2 & 5 & 0 & &  3 & 1 & 4\\
5 & 2 & 3 & &  0 & 4 & 1
\end{array}
\right]
\end{equation}
now with $R= 6849$ and
\begin{equation}
\label{isda3}
\mathbf{L}_{6D}^{\mathcal{R}}=\left[
\begin{array}
[c]{ccccccc}%
0 & 1 & 2 & &  3 & 4 & 5\\
1 & 0 & 3 & & 2 & 5 & 4\\
2 & 3 & 4 & & 5 & 0 & 1\\
   &    &   &  &   &    &   \\
3 & 2 & 5 & & 4 & 1 & 0\\
4 & 5 & 0 & & 1 & 2 & 3\\
5 & 4 & 1 & & 0 & 3 & 2
\end{array}
\right].
\end{equation}with $R=4689$. The dispersed $\mathbf{L}_{6D}^{\mathcal{R}}$ has the same set of eigenvalues as the dispersed $\mathbf{L}_{6D}$ and $\mathbf{L}_{6A}^{\mathcal{R}}$ has the same eigenvalues as $\mathbf{L}_{6A}$.  Note also the degeneracy of the singular values in both of the aggregated and dispersed compound Latin squares. In general, reverse compound squares share the same eigenvalues, singular values and rank as (\ref{cmpA}) and (\ref{cmpD}), which we will show rigorously in the following section. Thus, the compounding procedure for mixed order $m \neq n$ produces four unique squares, which are paired in terms of eigenvalues and singular values. For the degenerate order case, $m=n$, only one unique pair of squares can be generated.

While several authors, e.g. Derksen et al (\cite{derksen1}), Ollerenshaw
(\cite{oller2}), studying magic squares have effectively used the $k=2$
version of our associated type (\ref{cmpA}), we are not aware of any Kronecker
style studies of the distributed type square in (\ref{cmpD}).

Note that while we have derived the compounding formulae with matrices in mind, the Kronecker product is general to arbitrary tensors. This generality allows us to use the same formalism to describe the construction of Latin and magic cubes and hypercubes, and scales generally to arrays of arbitrary dimensionalty. We illustrate the compounding procedure for the simplest Latin cube, $N=2$, in Appendix \ref{appCube}.

\section{Generalization using GAPDA Arrays}
\label{sec:compoundingGAPDA}

While ISDA arrays allow a particular set of compound arrays to be generated, there is a larger class that can be formed without the requirement that the seed squares be ISDA themselves. In this case, we want to use arrays that are comprised of general arithmetic progressions, arranged to be doubly affine. We refer to these as GAPDA arrays. GAPDA arrays differ from ISDA arrays in that gaps can be present in the elements and the full cover of integers from $0..(n^k-1)$ is not required. However, to produce a compound magic array, the GAPDA seeds must be chosen such that the compound result produces an ISDA array with full coverage of the integer sequence and does not display any gaps. We do not discuss how to choose the individual GAPDA arrays which lead to full-cover products themselves. Instead, we will simply provide the details of the construction method and the properties that the compound ISDA array has, along with some examples of GAPDA seed pairs that generate ISDA arrays with full integer cover.

To generalize the compounding method, let us now assume that the matrices $A_m$ and $B_n$ are GAPDA arrays, and let us return to the Latin square example from Section \ref{sec:compoundingISDA}. For both the aggregated and dispersed compound squares, the elements in the resulting compound array are $\{0,1,2,3\}$. In this example we are looking for a compound square that has a tiling structure like $\mathbf{L}_2$ using the ISDA approach which has identical Latin seeds of $\mathbf{L}_2$. This Latin square has a full set of elements $\{0,1\}$ (ie, it is an ISDA square) and therefore must be multiplied by the coefficient $n$ in the compounding procedure. However, this need not be the case. We could equally well use an array with the same structure as $\mathbf{L}_2$, \textit{but with individual elements of $\{0,2\}$}. This is a GAPDA array since the sequence skips over $1$ entirely (See \ref{sec:ex3} for an explicit example). In this case, the multiplicative coefficients that were previously required to produce full cover compound squares (the $n^k$ and $m^k$ factors in equations \ref{cmpA} and \ref{cmpD}), are no longer explicitly needed. Thus, compounding requires we form the product
\begin{equation}
\text{ }
\mathbf{C}_{G}=
\left[
\begin{array}
[c]{ccc}
a_{11}\mathbf{E}_{n} + \mathbf{B}_{n}& ... & a_{1m}\mathbf{E}_{n}  + \mathbf{B}_{n} \\
\vdots & \ddots & \vdots \\
a_{m1}\mathbf{E}_{n}  + \mathbf{B}_{n} & ... & a_{mm}\mathbf{E}_{n}  + \mathbf{B}_{n}
\end{array}
\right].
\label{GAPDA}
\end{equation}
This form tiles $\mathbf{B}_{n}$ and individually augments each block by the elements of $\mathbf{A}_{m}$. In this case, the constant coefficients $n^k$ ($m^k$) are absorbed into the elements of the seed squares and no longer explicitly appear. (In the ISDA approach we have four seed squares $\{n^kA_m,B_n\}$ and $\{A_m,m^kB_n\}$). The distinction is subtle, but it significantly extends the generality of the compounding procedure, as it eliminates the need for full cover of the integers in the seed squares. This means we are free to consider other possible seeds which do not individually have full cover (i.e., GAPDA seeds), but which combine together to form fully covered ISDA compound arrays. Thus, the compounding form for GAPDA arrays appears even simpler than the ISDA case. Using one GAPDA seed pair $(\mathbf{A}_m,\mathbf{B}_n)$:
\begin{equation}
\mathbf{C}_{GA}= \left( \mathbf{A}_{m} \otimes \mathbf{E}_{n} \right)+ \left( \mathbf{E}_{m} \otimes \mathbf{B}_{n} \right).
\label{cmpGAP}
\end{equation}
The dispersed square, in contrast, always exhibits elements whose arrangement
cannot be restored by permutation to the aggregated array, and will come from a second distinct pair $(\mathbf{A}_m^\prime, \mathbf{B}_n^\prime)$:
\begin{equation}
\mathbf{C}_{GD}= \left( \mathbf{A^\prime}_{m} \otimes \mathbf{E}_{n} \right)+ \left( \mathbf{E}_{m} \otimes \mathbf{B^\prime}_{n} \right).
\label{cmpGAP}
\end{equation}
Reversing the order of the Kronecker products for each of the pair will produce perfectly shuffled squares. The total number of squares that can be formed using GAPDA squares varies by order. We show the number of squares that can be formed given GAP sequences in Table \ref{GAPTable}.

\begin{table}[htp]
\tiny
\linespread{1.0}
\center
\begin{tabular}{lllllll} \toprule
\scriptsize $\mathbf{N}$ & \scriptsize $\mathbf{C}_{mn}$ &\scriptsize $\mathbf{G}_{0}$ & \scriptsize $\mathbf{G}_{1}$ & \scriptsize $\mathbf{G}_{2}$ & \scriptsize $\mathbf{G}_{3}$ & \scriptsize $\mathbf{G}_{4}$
\\
\midrule
\\
\textbf{4} & \emph{2} & {$\{0,1\}$} &  &  &  &
\\
 &  & $\{0,2\}$ &  &  &  &
\\
\\
\textbf{6} & \emph{4} & $\{0,1,2\}$ & $\{0,2,4\}$ &  &  &
\\
 &  & $\{0,3\}$ & $\{0,1\}$ &  &  &
\\
\\
\textbf{8} & \emph{4} & $\{0,1,2,3\}$ & $\{0,2,4,6\}$ &  &  &
\\
 &  & $\{0,4\}$ & $\{0,1\}$ &  &  &
\\
\\
\textbf{9} & \emph{2} & $\{0,1,2\}$ &  &  &  &
\\
 &  & $\{0,3,6\}$ &  &  &  &
\\
\\
\textbf{10} & \emph{4} & $\{0,1,2,3,4\}$ & $\{0,2,4,6,8\}$ &  &  &
\\
 &  & $\{0,5\}$ & $\{0,1\}$ &  &  &
\\
\\
\textbf{12} & \emph{8} & $\{0,1,2,3,4,5\}$ & $\{0,2,4,6,8,10\}$ & $\{0,1,2,3\}$ & $\{0,3,6,9\}$ &
\\
 &  & $\{0,6\}$ & $\{0,1\}$ & $\{0,4,8\}$ & $\{0,1,2\}$ &
\\
\\
\textbf{14} & \emph{4} & $\{0,1,2,3,4,5,6$\} & $\{0,2,4,6,8,10,12\}$ &  &  &
\\
 &  & $\{0,7\}$ & $\{0,1\}$ &  &  &
\\
\\
\textbf{15} & \emph{4} & $\{0,1,2,3,4\}$ & $\{0,3,6,9,12\}$ &  &  &
\\
 &  & $\{0,5,10\}$ & $\{0,1,2\}$ &  &  &
\\
\\
\textbf{16} & \emph{10} & $\{0,1,2,3,4,5,6,7\}$ & $\{0,2,4,6,8,10,12,14\}$ & $\{0,4,8,12\}$ & $\{0,2,8,10\}$ & $\{0,1,8,9\}$
\\
 &  & $\{0,8\}$ & $\{0,1\}$ & $\{0,1,2,3\}$ & $\{0,1,4,5\}$ & $\{0,2,4,6\}$
\\
\bottomrule
\end{tabular}

\caption{GAPs for Compound Latin Squares of Small Composite Order.\smallskip \\
The generalized arithmetic progressions $\mathbf{G}_i$ provide complete cover for the composite order $N=mn$ squares producing a total of $\mathbf{C}_{mn}$ compounds for each pair of Latin squares used in their creation.}
\label{GAPTable}
\end{table}
\normalsize

To demonstrate the production of magic squares using GAPDAs, consider the two GAPDA arrays arranged in Lo Shu fashion,
\begin{equation}
\mathbf{K}_{0}=\left[
\begin{array}
[c]{ccc}%
9 & 20 & 1\\
2 & 10 & 18 \\
19 & 0 & 11 \\
\end{array}
\right]
\label{GAPDA1}
\end{equation}
and
\begin{equation}
\mathbf{K}_{1}=\left[
\begin{array}
[c]{ccc}%
27 & 60 & 3 \\
6 & 30 & 54 \\
57 & 0 & 33 \\
\end{array}
\right].
\label{GAPDA2}
\end{equation}
When these two squares are used as GAPDA seeds for the compounding formula given in equation (\ref{cmpGAP}), the result is the Koo-Soo-Ryak square (\ref{KSR}), given in \cite{KSR}. In this instance we can only form a single pair since $K_1$ is a simple multiple of $K_0$. By reversing the order of the Kronecker products we are still able to produce a second new square which is a perfect two dimensional shuffle of the original.

The utility of the GAPDA approach can be demonstrated by producing an additional pair of compound squares which does not seem to be present in the early historical record prior to the 20$^{th}$ century. If we consider the two compound arrays structured as the Lo Shu from the GAPs $\{0,3,6,9,12,15,18,21,24\}$ and $\{0,1,2,27,28,29,54, 55,56\}$ we arrive at the new compound square $\mathbf{C}_{B}$ and its shuffle first seen in Frierson, one of the pair being:
\begin{equation}
\mathbf{C}_{B}=
\left[
\begin{array}
[c]{ccccccccccc}
36 & 65 & 10 & & 51 & 80 & 25 & & 30 & 59 &  4 \\
11 & 37 & 63 & & 26 & 52 & 78 & &  5 & 31 & 57 \\
64 &  9 & 38 & & 79 & 24 & 53 & & 58 &  3 & 32 \\
&&&&&&&&&&\\
33 & 62 &  7 & & 39 & 68 & 13 & & 45 & 74 & 19 \\
 8 & 34 & 60 & & 14 & 40 & 66 & & 20 & 46 & 72 \\
61 &  6 & 35 & & 67 & 12 & 41 & & 73 & 18 & 47 \\
&&&&&&&&&&\\
48 & 77 & 22 & & 27 & 56 &  1 & & 42 & 71 & 16 \\
23 & 49 & 75 & & 2  & 28 & 54 & & 17 & 43 & 69 \\
76 & 21 & 50 & & 55 &  0 & 29 & & 70 & 15 & 44
\end{array}
\right]
\end{equation}
This now accounts for all six of the basic possible compound magic squares of order nine in standard form discussed in Frierson \cite{Frierson} (with R values of $[1,301,165,856]$, $[797,281,056]$ and $[842,630,688]$). Notice that we have three pairs of seed squares of the same order for the magic square in contrast to just one pair for the Latin square, due to the fact that the elements in the final compound square are $0\ldots (n^2-1)$ rather than $0\ldots (n-1)$.

\section{Compound matrix eigenspectra}
\label{sec:eigenspectra}

Let us return to the ISDA case to discuss the compound matrix eigenspectra. Consider the definition of the Kronecker product for matrices in \cite{Conway}. Denote the set of $n$ eigenvectors of $\mathbf{B}_{n}$ as $\mathbf{u}_{t}$ where $1 \leq t \leq n$. Each of these eigenvectors has elements $\mathbf{u}_{t}=(u_{t1}, u_{t2}, ..., u_{tn})^T$. The corresponding eigenvalues $\mu_{t}$ satisfy
\begin{equation}
\mathbf{B}_{n} \mathbf{u}_{t} = \mu_{t} \mathbf{u}_{t},
\label{eigBmu}
\end{equation}
by definition. Note that for the time being, we make no assumptions about the structure or properties of $\mathbf{B}_{n}$, or the values of the individual $\mu_{t}$. The matrix $\mathbf{E}_{m}$ has only one non-zero eigenvalue, $m$, with eigenvector $\mathbf{e}_{m}$, a vector with each element equal to $1$. Due to these relationships we write the eigenvectors of
$\mathbf{\beta}_{n}$ in terms of the Kronecker product:
\begin{equation}
\mathbf{\Psi}_{t}=\mathbf{e}_{m} \otimes \mathbf{u}_{t}=
\left[
\begin{array}
[c]{c}
\mathbf{u}_{t}\\
\vdots\\
\mathbf{u}_{t}
\end{array}
\right].
\label{psiB}
\end{equation}
Now consider the relationship between the eigenvalues and eigenvectors of $\mathbf{\beta}_{n}$. Due to the tiled structure of this matrix we see that
\begin{equation}
\mathbf{\beta}_n \mathbf{\Psi}_{t}=\Lambda_{t}\mathbf{\Psi}_{t},
\label{BBBB}
\end{equation}
where the eigenvalues of this tiled matrix are
\begin{equation}
\Lambda_{t} = m \mu_{t},
\label{LambdaEig}
\end{equation}
as seen from the tiled structure of $\mathbf{\beta}_n$ in (\ref{beta}). Note that $\mathbf{\beta}_{n}$ is singular since it has at most only $n$ eigenvalues as determined by the spectrum of $\mathbf{B}_{n}$.

Now let us consider $\mathbf{\alpha}_m$, given by (\ref{alpha}). The relationship between the eigenvalues and eigenvectors of $\mathbf{A}_{m}$ is defined in analogy with the notation defined above with $1 \leq s \leq m$, such that
\begin{equation}
\mathbf{A}_{m}\mathbf{v}_{s}=\nu_{s}\mathbf{v}_{s},
\label{eigAlambda}
\end{equation}
where the eigenvalues of $\mathbf{A}_{m}$ are $\nu_{s}$ and the
eigenvectors have elements $\mathbf{v}_{s}=(v_{s1},v_{s2},...,v_{sm})^T$. The matrix $\mathbf{\alpha}_{m}$ has a set of eigenvectors that can be expressed in terms of the Kronecker product
\begin{equation}
\mathbf{\Phi}_{s}=\mathbf{v}_{s} \otimes \mathbf{e}_{n}\text{ }=
\left[
\begin{array}
[c]{c}
v_{s1} \mathbf{e}_{n}\\
\vdots\\
v_{sm} \mathbf{e}_{n}
\end{array}
\right],
\label{phi}
\end{equation}
such that the corresponding eigenvalues of $\mathbf{\alpha}_{m}$ are
\begin{equation}
\mathbf{\alpha}_{m} \mathbf{\Phi}_{s}=\Delta_{s} \mathbf{\Phi}_{s},
\label{aaaa}
\end{equation}
with
\begin{equation}
\Delta_{s}=n \nu_s.
\label{DeltaEig}
\end{equation}
In analogy with $\mathbf{\beta}_n$, $\mathbf{\alpha}_{m}$ is also singular in general, and has at most $m$ non-trivial eigenvalues.

The GAPDA matrices display analogous behaviour, only with the factors of $n^k$ ($m^k$) eliminated.

\subsection{The eigenproperties of semi-magic matrices}
\label{secEigSMM}

To describe the eigenproperties of ISDA compound squares we require merely that $\mathbf{A}_{m}$ and $\mathbf{B}_{n}$ be doubly-affine arrays. This condition introduces structure into the spectra of the seed matrices and significantly simplifies the general results presented above. Non-negative semi-magic matrices have a dominant eigenvalue, which we denote by $\nu_{1}$ and $\mu_1$, respectively. These dominant eigenvalues set an upper bound on the eigenvalue spectrum for each constituent square. The remaining eigenvalues are the ``non-dominant'' eigenvalues which we label with $2 \leq s \leq m$ and $2\leq t \leq n$ for $\mathbf{A}_m$ and $\mathbf{B}_n$. We then have $\left\vert \nu_{1}\right\vert > \left\vert \nu_{t}\right\vert$ and $\left\vert \mu_{1}\right\vert > \left\vert \mu_{s}\right\vert$. This is a consequence of Perron's theorem \cite{Ortega}. For magic and Latin matrices, the dominant eigenvalues are equal to the line sums (\ref{magicSum}) of the constituent squares:
\begin{equation}
\nu_{1}=S_{m,k}
\label{dominantLA}
\end{equation}
\begin{equation}
\mu_{1}=S_{n,k}
\label{dominantLB}
\end{equation}
and by the doubly-affine property, the eigenvectors $\mathbf{v}_{1}
=\mathbf{e}_{m}$ and $\mathbf{u}_{1}=\mathbf{e}_{n}$ associated with the
dominant eigenvalues have a $1$ at each position. This is both a left and right eigenvector (i.e. for doubly-affine matrices and their transposes), so the
principle of bi-orthogonality \cite{Ortega} ensures the other, non-dominant
eigenvectors on both the left and the right will be orthogonal to
$\mathbf{e}_{m}$ and $\mathbf{e}_{n}$. This means that the non-dominant
eigenvectors of doubly-affine matrices have elements which sum to zero in general.

The semi-magic properties have consequences for the eigenproperties of the Kronecker products $\mathbf{\beta}_n$ and $\mathbf{\alpha}_m$ derived earlier. The dominant eigenvalues of $\left(\ref{aaaa}\right)  $ and $\left(  \ref{BBBB}\right) $ are the line sums  $\Delta_1=nS_{m,k}$ and $\Lambda_1=mS_{n,k}$. The eigenvectors $\mathbf{\Phi}_1$ and $\mathbf{\Psi}_1$ associated with these dominant eigenvalues are $\mathbf{e}_{nm}$ for both squares by the arguments above.

Now let us return to the remaining non-dominant eigenvectors of $\mathbf{\alpha}_{m}$ and $\mathbf{\beta}_{n}$. The non-dominant eigenvectors are $\mathbf{\Phi}_s$ and $\mathbf{\Psi}_t$ are defined in $\left(\ref{phi}\right)  $ and $\left(  \ref{psiB}\right)  $. Let us consider the effect that $\mathbf{\alpha}_{m}$ has on the non-dominant eigenvectors of $\mathbf{\beta}_{n}$:
\begin{equation}
\mathbf{\alpha}_{m}\mathbf{\Psi}_{t}=\mathbf{0},
\label{adam12}
\end{equation}
where $\mathbf{0}$\ is the $nm$ element null vector. From the the definition, we see a similar effect when $\mathbf{\beta}_{n}$ acts on the eigenvectors of $\mathbf{\alpha}_{m}$:
\begin{equation}
\mathbf{\beta}_{n} \mathbf{\Phi}_{s}=\mathbf{0}.
\label{adam13}
\end{equation}
These relationships hold due to the form of the eigenvectors defined in $\left(\ref{phi}\right)  $ and $\left(  \ref{psiB}\right)  $ and the tiled structures of $\mathbf{\beta}_n$ and $\mathbf{\alpha}_m$, and are a consequence of the principle of bi-orthogonality since the sum of the elements of the non-dominant eigenvectors must vanish, i.e. $\sum_{j=1}^{n} u_{tj}=0$  for all $t \neq 1$ and $\sum_{j=1}^{m} v_{sj}=0$ for all $s \neq 1$. We can therefore conclude that the non-dominant eigenvectors of $\mathbf{\beta}_{n}$ occupy the null space of $\mathbf{\alpha}_{m}$ and vice-versa.

The dominant eigenvector of both $\boldsymbol{\alpha}_{m}$ and
$\boldsymbol{\beta}_{n}$ is $\mathbf{e}_{nm}$, and must also be an eigenvector
of the compound square. The dominant eigenvalue of $\mathbf{C}_{A}$ is
$\lambda_{1}=\Lambda_{1}+n^{k}\Delta_{1}$. However, the relationships given in
(\ref{LambdaEig}), (\ref{DeltaEig}), (\ref{dominantLB}) and (\ref{dominantLA})
allow us to simplify this expression considerably:
\begin{equation}
\lambda_{1}=S_{mn,k},
\end{equation}
which has eigenvector $\mathbf{e}_{nm}$ as expected for an order $mn$
doubly-affine matrix. The remaining $n+m-1$ possibly non-vanishing,
non-dominant eigenvalues of the aggregated compounded square are $m\mu_{t}$
and $n^{k+1}\nu_{s}$, with eigenvectors given by $\left(  \ref{psiB}\right)  $
and $\left(  \ref{phi}\right)  $, respectively by the previous arguments
presented in this section. A similar analysis can be performed for the
dispersed compound square $\mathbf{C}_{D}$ given by (\ref{cmpD}), which yields
the remaining potentially non-vanishing non-dominant eigenvalues $m^{k+1}
\mu_{t}$ and $n\nu_{s}$. Note that these arguments require the non-degeneracy
of the eigenvalues. In the case of degeneracy, the eigenvalue formulae will
still be applicable though the description of the eigenvectors will be more
complicated, and this case is not addressed here. A similar issue arises when
discussing the singular values in Section \ref{sec:SVD}.

The reverse aggregated compound square $\mathbf{C}_A^{\mathcal{R}}$ (\ref{revCmpA}) and reverse dispersed square $\mathbf{C}_D^{\mathcal{R}}$ (\ref{revCmpD}) have similar eigenvalue properties as the aggregated and dispersed squares $\mathbf{C}_A$ and $\mathbf{C}_D$ which can be shown by similar arguments to those used above.

\section{Singular Value Decomposition}
\label{sec:SVD}

The singular value decomposition (SVD; \cite{GR}, \cite{GolVan}) expresses a matrix $\mathbf{M}$ as a product $\mathbf{M}=\mathbf{U \Sigma V}^T$, where $\mathbf{U}$ and $\mathbf{V}$ are orthogonal matrices and $\mathbf{\Sigma}$ is diagonal, with non-zero elements that define the singular values. The columns of $\mathbf{U}$ and $\mathbf{V}$ are the left and right eigenvectors of the products $\mathbf{MM}^T$ and $\mathbf{M}^T\mathbf{M}$, respectively. We refer to these matrix products interchangably as the Gramian matrix. In general the eigenvalues of both of these products are identical, and are defined as the squares of the singular values. Cameron, Rogers and Loly
(CRL) \cite{ian2},\cite{CRL} used the singular values to calculate the Shannon
entropy (see Newton and De Salvo \cite{NDS}) of doubly-affine squares to
understand these matrices. The SV ``clans'' CRL introduced express the
significance and utility of the SVs over eigenvalue sets. The Shannon entropy,
recently used by \cite{NDS} in the context of Latin squares and Sudokus and
applied to magic squares in \cite{CRL}, can then be evaluated using the SV
formulae given in the present work. CRL's calculations show that compound
squares usually have lower entropy than most other squares of the same order.

Just as we found a relationship between the eigenvalues of the compound and seed squares, a similar relationship exists between the singular values. In analogy with the arguments used to derive the eigenvector properties, we can perform a similar calculation using the singular value basis vectors. Let us define the set of non-dominant singular values of $\mathbf{A}_m$ as $\xi_s$, and $\mathbf{B}_n$ as $\eta_t$. The indices $s$ and $t$ are defined as in the previous section. We define the singular values of the compound square $\sigma_i$. Consider the Gramian of the aggregated compound square:
\begin{equation}
\mathbf{C}_{A}^{T} \mathbf{C}_{A}=\left(  n^{k}\mathbf{\alpha
}_{m}+\mathbf{\beta}_{n}\right)  ^{T}\left(  n^{k}\mathbf{\alpha}
_{m}+\mathbf{\beta}_{n}\right)
\end{equation}
such that
\begin{equation}
\mathbf{C}_{A}^{T} \mathbf{C}_{A}=n^{2k}\mathbf{\alpha}_{m}^{T}\mathbf{\alpha}_{m}+n^{k}\mathbf{\beta}_{n}^{T}\mathbf{\alpha}_{m}
+n^{k}\mathbf{\alpha}_{m}^{T}\mathbf{\beta}_{n}+\mathbf{\beta}_{n}^{T}\mathbf{\beta}_{n}. \label{Gram}
\end{equation}
The mixed terms $\mathbf{\beta}_n^{T}\mathbf{\alpha}_m$, $\mathbf{\alpha
}_m^{T}\mathbf{\beta}_n$ are each others transpose and only contribute to the compound line sum (and largest) singular value with singular vector $\mathbf{e}_{mn}$:
\begin{equation}
\sigma_{1}=\left[ (mS_{n,k})^2+ 2n^k(mnS_{m,k}S_{n,k}) + n^{2k}(nS_{m,k})^2 \right]^{\frac{1}{2}}=S_{C,k},
\end{equation}
The mixed terms in $(\ref{Gram})$ vanish when operating on non-dominant singular basis vectors due to the relationships in $(\ref{adam12})$ and $(\ref{adam13})$. The $m+n-1$ remaining possibly non-zero singular values are $n^{k+1} \xi_s$ and $m \eta_t$. The non-mixed terms in (\ref{Gram}) are the Gramians of $\mathbf{\alpha}$ and $\mathbf{\beta}$ respectively.

In general, these expressions for the singular values apply for both magic and Latin squares. However, we have observed degeneracy in the singular value spectra of Latin squares ($k=1$). The presence of the degeneracy of the singular values for some squares complicates the analytical description of the SVD basis. Despite this, the arguments that we have presented above hold for the singular values even in the case of degeneracies. Just as  the eigenvalue case, the GAPDA case follows by similar arguments.

We summarize the results in Tables \ref{table1} and \ref{table2}, which list the ISDA and GAPDA results for the eigenvalues, eigenvectors and singular values.

\begin{table}[htp]
\small
\centering
\bgroup
\def\arraystretch{1.3}
\begin{tabular}{llllll}
\toprule
Square & Order & Formula & Eigenvalues & Eigenvectors & Singular Values \\
\midrule
$\mathbf{A}_m$ & m & Input & $S_{m,k}$ & $\mathbf{e}_m$ & $S_{m,k}$   \\
&     &          & $\nu_s$ & $\mathbf{v}_s$  & $\xi_s$  \\
$\mathbf{B}_n$  & n  & Input & $S_{n,k}$ & $\mathbf{e}_n$  & $S_{n,k}$ \\
&     &          & $\mu_t$ &  $\mathbf{u}_t$  & $\eta_t$ \\
\midrule
&     &          & $S_{mn,k}$ & $\mathbf{e}_{mn}$ & $S_{mn,k}$ \\

$\mathbf{C}_A$ & mn & $\beta_n + n^k \alpha_m$ & $m\mu_t$ & $\mathbf{e}_m \otimes \mathbf{u}_t$ & $m \eta_t$ \\
& &  & $n^{k+1}\nu_s$ & $\mathbf{v}_s \otimes \mathbf{e}_n$ & $n^{k+1} \xi_s$ \\
\midrule

&      &         & $S_{mn,k}$ & $\mathbf{e}_{mn}$ & $S_{mn,k}$ \\

$\mathbf{C}_D$ & mn & $\alpha_m + m^k \beta_n$ & $m^{k+1}\mu_t$ & $\mathbf{u}_t \otimes \mathbf{e}_m$ & $m^{k+1} \eta_t$ \\
& & & $n \nu_s$ & $\mathbf{e}_n \otimes \mathbf{v}_s$ & $n \xi_s$ \\
\bottomrule
\end{tabular}
\egroup
\caption{Summary of Compounding Results for ISDA\smallskip \\
Note that some of the non-dominant eigenvalues $\nu_{s}$, $\mu_{t}$ of the seed squares may vanish (see Example \ref{sec:ex4}). This is also true for the singular values, so that the rank of the compound square is further reduced. The necessarily zero eigenvalues and singular values of $\mathbf{C}_A$ and $\mathbf{C}_D$ are not displayed in the table. These results also hold for the reverse squares, $\mathbf{C}_A^{\mathcal{R}}$ and $\mathbf{C}_D^{\mathcal{R}}$. As previously defined, we take the indices $2 \leq s \leq m$ and $2 \leq t \leq n$.\\ \bigskip}
\label{table1}

\small
\bgroup

\def\arraystretch{1.3}
\begin{tabular}{llllll}
\toprule
Square & Order & Formula & Eigenvalues & Eigenvectors & Singular Values \\
\midrule
$\mathbf{A}_\text{m}^\prime$ & m & Input & $S_{m,k}^\prime$ & $\mathbf{e}_\text{m}$ & $S_{m,k}^\prime$   \\

&     &          & $\nu_\text{s}^\prime$ & $\mathbf{v}_\text{s}$  & $\xi_\text{s}^\prime$  \\

$\mathbf{B}_\text{n}^\prime$  & n  & Input & $S_{n,k}^\prime$ & $\mathbf{e}_{n}$  & $S_{n,k}^\prime$ \\
&     &          & $\mu_{t}^\prime$ &  $\mathbf{u}_{t}$  & $\eta_{t}^\prime$ \\
\midrule

&     &          & $S_{mn,k}$ & $\mathbf{e}_{mn}$ & $S_{mn,k}$ \\

$\mathbf{C}_{mn}$ & mn & $\beta_{n}^\prime + \alpha_{m}^\prime$ & $m\mu_{t}^\prime$ & $\mathbf{e}_{m} \otimes \mathbf{u}_{t}$ & $m \eta_{t}^\prime$ \\

& &  & $n\nu_{s}\prime$ & $\mathbf{v}_{s} \otimes \mathbf{e}_{n}$ & $n \xi_{s}\prime$ \\
\bottomrule

\end{tabular}
\egroup
\caption{Summary of Compounding Results for GAPDA\smallskip \\
Latin compounds of the form $m\times n$ where $m$
and $n$ are distinct have two pairs of GAP's available: one with $m$ elements, $A_{0}=\{0,1,\ldots,m-1\}$ and one  with $n$ elements, $B_{0}=\{0,m,\ldots,(n-1)m\}$
 and similarly $A_{1}=\{0,n,\ldots,(m-1)n\}$, $B_{1}=\{0,1,\ldots,n\}$.
This will always allow us to produce four compound squares and the distinction between aggregated and dispersed becomes arbitrary. Note that $n\nu_0+m\mu_0=S_{mn,k}$, where $\nu_0$ and $\mu_0$ are the dominant eigenvalues (row/column sums) of the individual seed squares.}
\label{table2}
\end{table}
\newpage
\section{Conclusions}
We have shown that it is possible to produce large order magic and Latin arrays in arbitrary dimensions through a compounding process that is described by the addition of kronecker products. For the two dimensional case of magic and Latin squares, the eigenvalues and eigenvectors of the seed squares are directly related to the eigenproperties of the compound squares in a simple way, and a similar relationship exists for the SVD of the compound square. There are four unique compound squares that can be formed for each pair of distinct order constituent ISDA squares. The aggregated $\mathbf{C}_A$ and reverse-aggregated $\mathbf{C}_A^{\mathcal{R}}$ squares group consecutive integers and have identical spectral and singular properties. The dispersed $\mathbf{C}_D$ and reverse-dispersed $\mathbf{C}_D^{\mathcal{R}}$ squares disperse consecutive integers throughout the square and share similar spectral and singular properties. The properties of compounding from the GAPDA matrices follow directly from the ISDA approach. We summarize the principal spectral and singular properties of compound squares in Tables \ref{table1} and \ref{table2}. Using the compounding formulae for both ISDA and GAPDA seed squares, we recovered the six order $9$ compound magic squares discussed by Frierson \cite{Frierson}.

We have found that Latin squares can possess degenerate singular values. In general this degeneracy complicates the analysis of the compound singular vectors. Therefore we leave an in-depth analysis of the degenerate singular vectors as an open problem.

\begin{appendices}
\section{Appendix: Example Problems}
In this appendix, we present several examples of compounding that highlight the eigenvalues, singular values and rank calculations. Note that we usually list \textit{squared} singular values $\sigma_i^2$.

\subsection{Spectral properties of a Latin compound square}
Let us return to the rank $r=2$ Latin square example from section \ref{sec1}
\begin{equation}
\mathbf{L}_2=\left[
\begin{array}
[c]{cc}
0 & 1\\
1 & 0
\end{array}
\right].
\end{equation}
The dominant eigenvalue is $\mu_1=1$ with eigenvector $\mathbf{u}_1=\mathbf{e}_2=(1,1)^T$ due to the semi-magic property, and non-dominant $\mu_2=-1$ with $\mathbf{u}_2=(-1,1)^T$. Note that the sum of the eigenvalues is zero, since the trace of $\mathbf{L}_2$ vanishes. The singular values of this square are degenerate, with $\sigma_1^2=\sigma_2^2=1$, and rank $r=2$. Compounding this square with itself sets $n=m=2$, $\nu_i=\mu_i$, $\xi_i=\eta_i$, and since the square is latin, $k=1$.
While the corresponding rank 1 augmentation matrix, $\mathbf{E}_{2}$, is
\begin{equation}
\mathbf{E}_{2}=\left[
\begin{array}
[c]{cc}
1 & 1\\
1 & 1
\end{array}
\right].
\end{equation}

The dispersed compound square of rank $r=3\textit{(=2+2-1)}$ is
\begin{equation}
\mathbf{C}_{D}=\left[
\begin{array}
[c]{cccc}%
0 & 2 & 1 & 3\\
2 & 0 & 3 & 1\\
1 & 3 & 0 & 2\\
3 & 1 & 2 & 0
\end{array}
\right].
\end{equation}
Since in this case $\mathbf{A}_m=\mathbf{B}_n$, the reverse aggregated square and the dispersed square are the same, $\mathbf{C}_A^{\mathcal{R}}=\mathbf{C}_D$ and  $\mathbf{C}_D^{\mathcal{R}}=\mathbf{C}_A$. Note that both of these squares have the same rank and are their own transposes.

\subsection{Mixed Order Compounding  - GAPDA Approach}
\label{sec:ex3}

Consider the latin squares:
\begin{equation}
\mathbf{L}_2=\left[
\begin{array}
[c]{ccc}%
0 & 1 \\
1 & 0
\end{array}
\right]
\hspace{5mm}\&\hspace{5mm}
\mathbf{L}_3=\left[
\begin{array}
[c]{ccc}
0 & 1 & 2\\
1 & 2 & 0\\
2 & 0 & 1
\end{array}
\right].
\end{equation}
with orders $m=2$ and $n=3$ respectively. Suppose we wish to produce a compound Latin square of order $m\times n$ that reflects the structure of both squares but has complete cover, ie this square will have elements $\{0,1,2,3,4,5\}$. All possible combinations of the elements selected from each of the two generalized arithmetic progressions (GAP's) $\{0,1,2\}$ and $\{0,3\}$ will generate the required cover set. (We see that the $0$ from the second set added to all the elements of the first will generate the first half of the cover set $\{0,1,2\}$ and the second element $3$ added to all the elements of the first set produces the second half of the cover set, $\{3,4,5\}$). Arranging the elements of the two GAP’s to reflect the requisite Latin square structures of $L_3$ and $L_2$ yields the two corresponding GAPDA arrays:
\begin{equation}
\label{G01}
\mathbf{G}_{0}=\left[
\begin{array}
[c]{ccc}%
0 & 1 & 2 \\
1 & 2 & 0 \\
2 & 0 & 1
\end{array}
\right]
\hspace{5mm}\&\hspace{5mm}
\mathbf{G}_{1}=\left[
\begin{array}
[c]{ccc}%
0 & 3 \\
3 & 0
\end{array}
\right]
\end{equation}
We find the resulting aggregated compound Latin square tiled as in (\ref{isda0}).

We could equally well use a different pair of GAP's with elements $\{0,2,4\}$ and $\{0,1\}$ to generate the cover set and use the corresponding doubly affine arrays:
\begin{equation}
\label{G23}
\mathbf{G}_{2}=\left[
\begin{array}
[c]{ccc}
0 & 2 & 4\\
2 & 4 & 0\\
4 & 0 & 2
\end{array}
\right]
\hspace{5mm}
\&\hspace{5mm}
\mathbf{G}_{3}=\left[
\begin{array}
[c]{ccc}%
0 & 1 \\
1 & 0
\end{array}
\right]
\end{equation}
to yield the dispersed square as in (\ref{isda1}) which appears to be similar in structure to the previous compound square \textit{but has a different ordering of the individual elements}.

Simply reversing the order of operations with the previous two pair of GAPDA's (\ref{G01}) and (\ref{G23}):
\begin{equation}
\mathbf{G}_{1}=\left[
\begin{array}
[c]{ccc}%
0 & 3 \\
3 & 0
\end{array}
\right]
\hspace{5mm}
\&
\hspace{5mm}
\mathbf{G}_{0}=\left[
\begin{array}
[c]{ccc}%
0 & 1 & 2 \\
1 & 2 & 0 \\
2 & 0 & 1
\end{array}
\right]
\end{equation}\\
and
\begin{equation}
\mathbf{G}_{3}=\left[
\begin{array}
[c]{ccc}%
0 & 1 \\
1 & 0
\end{array}
\right]
\hspace{5mm}
\&
\hspace{5mm}
\mathbf{G}_{2}=\left[
\begin{array}
[c]{ccc}%
0 & 2 & 4 \\
2 & 4 & 0 \\
4 & 0 & 2
\end{array}
\right]
\end{equation}

we can form the reverse aggregated and dispersed pair, $\mathbf{C}_{A}^{\mathcal{R}}$ \& $\mathbf{C}_{D}^{\mathcal{R}}$ as in (\ref{isda2}) and (\ref{isda3}) which are \textit{perfect shuffles of their earlier counterparts and which exhibit a different tile structure based on the seeds $\mathbf{G}_{1}$ and $\mathbf{G}_{3}$}. We show the aggregated pair $\mathbf{C}_A$, $\mathbf{C}_A^{\mathcal{R}}$:
\begin{equation}
\label{duet1}
\left[
\begin{array}
[c]{ccccccc}%
\textcolor{blue}0 & 1 & 2 & & \textcolor{blue}3 & 4 & 5\\
1 & 2 & 0 & & 4 & 5 & 3\\
2 & 0 & 1 & & 5 & 3 & 4\\
  &   &   & &   &   &  \\
\textcolor{blue}3 & 4 & 5 & & \textcolor{red}0 & \textcolor{red}1 & \textcolor{red}2\\
4 & 5 & 3 & & \textcolor{red}1 & \textcolor{red}2 & \textcolor{red}0\\
5 & 3 & 4 & & \textcolor{red}2 & \textcolor{red}0 & \textcolor{red}1
\end{array}
\right]\rightleftharpoons
\left[
\begin{array}
[c]{cccccccc}
\textit{\textcolor{blue}0} & 3 & & \textit{\textcolor{blue}1} & 4 & & \textit{\textcolor{blue}2} & 5\\
3 & \textit{0} & & 4 & \textit{1} & & 5 & \textit{2}\\
  &   & &   &   & &   &  \\
\textit{\textcolor{blue}1} & 4 & & \textit{\textcolor{blue}2} & 5 & & \textit{\textcolor{blue}0} & 3\\
4 & \textit{1} & & 5 & \textit{2} & & 3 & \textit{0}\\
  &   & &   &   & &   &  \\
\textit{\textcolor{blue}2} & 5 & & \textit{\textcolor{red}0} & \textcolor{red}3 & & \textit{\textcolor{blue}1} & \textit{4}\\
5 & \textit{2} & & \textcolor{red}3 & \textit{\textcolor{red}0} & & 4 & \textit{1}
\end{array}
\right]
\end{equation}

and the dispersed pair $\mathbf{C}_D$, $\mathbf{C}_D^{\mathcal{R}}$:

\begin{equation}
\label{duet2}
\left[
\begin{array}
[c]{ccccccc}%
\textcolor{blue}0 & 2 & 4 & & \textcolor{blue}1 & 3 & 5\\
2 & 4 & 0 & & 3 & 5 & 1\\
4 & 0 & 2 & & 5 & 1 & 3\\
  &   &   & &   &   &  \\
\textcolor{blue}1 & 3 & 5 & & \textcolor{red}0 & \textcolor{red}2 & \textcolor{red}4\\
3 & 5 & 1 & & \textcolor{red}2 & \textcolor{red}4 & \textcolor{red}0\\
5 & 1 & 3 & & \textcolor{red}4 & \textcolor{red}0 & \textcolor{red}2
\end{array}
\right]
\rightleftharpoons
\left[
\begin{array}
[c]{cccccccc}%
\textit{\textcolor{blue}0} & 1 & & \textit{\textcolor{blue}2} & 3 & & \textit{\textcolor{blue}4} & 5\\
1 & \textit{0} & & 3 & \textit{2} & & 5 & \textit{4}\\
  &   & &   &   & &   &  \\
\textit{\textcolor{blue}2} & 3 & & \textit{\textcolor{blue}4} & 5 & & \textit{\textcolor{blue}0} & 1\\
3 & \textit{2} & & 5 & \textit{4} & & 1 & \textit{0}\\
  &   & &   &   & &   &  \\
\textit{\textcolor{blue}4} & 5 & & \textit{\textcolor{red}0} & \textcolor{red}1 & & \textit{\textcolor{blue}2} & 3\\
5 & \textit{4} & & \textcolor{red}1 & \textit{\textcolor{red}0} & & 3 & \textit{2}
\end{array}
\right].
\end{equation}
An inspection of the location of the italicized numerals in the reversed squares shows that row and column permutations will restore the previous upper and lower corner blocks and the ordering of the diagonal elements in the aggregated and dispersed squares. So we \textit{know} that the
dispersed $\mathbf{C}_D^{\mathcal{R}}$ will have the same set of singular values as the dispersed $\mathbf{C}_D$ and
$\mathbf{C}_A^{\mathcal{R}}$ to have the same singular values as $\mathbf{C}_{A}$. \\
The arrays in (\ref{duet1}) and (\ref{duet2}) nicely exhibit 2-dimensional perfect multiway shuffles, first a 2-way shuffle proceeding from left to right and then a 3-way shuffle proceeding from right to left. Consider the top row in the first square in (\ref{duet2}), which we had divided to reflect the tile structure of the square. With a deck of cards the 2-way perfect or riffle shuffle proceeds by splitting the deck into two halves or subdecks and alternately interleaving the cards from each half to produce the full shuffled deck. In our case the two decks represented by $[024\cdot 135]$ through the shuffling process produce the sequence $[012345]$, as can be seen in the top row of the second square. This operation affects the column ordering. Since we need a 2-dimensional shuffle, we next have to consider the shuffling operation on the rows by examining the columns. The first column $[024\cdot 135]$ through the shuffling process produces the sequence $[012345]$, as before. The net result of repeating the entire process for all the rows and columns generates the second square which we had divided to explicitly show the new tile structure. The 3-way shuffle proceeds in similar fashion, notice that the top row of the second square already conveniently divided $[01\cdot 23\cdot 45]$ into 3 decks goes to $[024\cdot 135]$ by selecting a corresponding element from each deck in turn.

\subsection{An order $12$ aggregated magic square}
\label{sec:ex4}
To produce an order $12$ magic square, we will compound the order $n=3$ Lo Shu $\mathbf{M}_3$ as given in (\ref{M3}) with an order $m=4$ magic square. Of the many possible choices that exist among the standard 880, which exclude rotations and reflections, we use:
\begin{equation}
\mathbf{M}_4=\left[
\begin{array}
[c]{cccc}
4 & 3 & 15 & 8 \\
10 & 13 & 1 & 6 \\
9 & 14 & 2 & 5 \\
7 & 0 & 12 & 11
\end{array}
\right].
\end{equation}
This square, which we find in Loly et al \cite{LCTS} has a single non-zero eigenvalue that corresponds to the line sum, $\nu_1=30$ \cite{powers}. Despite the single non-zero eigenvalue of $\mathbf{M}_4$, it is a rank $3$ square with singular values $\xi_1^2=\nu_1^2$, $\xi_2^2=320$ and $\xi_3^2=20$. The order $3$ magic square has eigenvalues $\mu_1=12$, $\mu_2=2i\sqrt{3}$ and $\mu_3=-2i \sqrt{3}$, and singular values $\eta_i^2=\{144,48,12\}$.

The aggregated compound square is
\small
\begin{equation}
\mathbf{C}_A=\left[
\begin{array}
[c]{cccccccccccccccc}
39& 44&
37& 30&
35& 28&
138& 143&
136& 75&
80& 73\\
38& 40&
42& 29&
31& 33&
137& 139&
141& 74&
76& 78\\
43& 36&
41& 34&
27& 32&
142& 135&
140& 79&
72& 77\\
93& 98&
91& 120&
125& 118&
12& 17&
10& 57&
62& 55\\
92& 94&
96& 119&
121& 123&
11& 13&
15& 56&
58& 60\\
97& 90&
95& 124&
117& 122&
16& 9&
14& 61&
54& 59\\
84& 89&
82& 129&
134& 127&
21& 26&
19& 48&
53& 46\\
83& 85&
87& 128&
130& 132&
20& 22&
24& 47&
49& 51\\
88& 81&
86& 133&
126& 131&
25& 18&
23& 52&
45& 50\\
66& 71&
64& 3&
8& 1&
111& 116&
109& 102&
107& 100\\
65& 67&
69& 2&
4& 6&
110& 112&
114& 101&
103& 105\\
70& 63&
68& 7&
0& 5&
115& 108&
113& 106&
99& 104
\end{array}
\right].
\end{equation}
\normalsize
This order $12$ square is rank $r=5\textit{(=3+3-1)}$ exactly as we would expect from Tables \ref{table1} and \ref{table2} where we see that $rank(\mathbf{C}) = rank(\mathbf{A}) + rank(\mathbf{B}) -1$, and has eigenvalues $\lambda_1=858$, $\lambda_2=4(2i\sqrt{6})$, $\lambda_3=-4(2i\sqrt{6})$. The singular values are $\sigma_i^2=\{858^2,9^2\times320,9^2\times20,4^2\times48,4^2\times12\}$.

\subsection{Compounding a Latin Cube}
\label{appCube}

Here we demonstrate the compounding process for an $N=2$ Latin cube. This is an ISDA array with $k=1$:
\begin{equation}
\textbf{L}_\text{cube}=
\begin{tabular}{ll}
$\fbox{$%
\begin{array}{rr}
0 & 1\\
1 & 0%
\end{array}%
$}$ & $\fbox{$%
\begin{array}{rr}
1 & 0 \\
0 & 1 %
\end{array}%
$}$ \\
\multicolumn{1}{c}{1$^{\text{st}}$ Layer} & \multicolumn{1}{c}{2$^{\text{nd}}$ Layer}%
\end{tabular}
\label{latinCubeSeed}
\end{equation}
Using eq. \ref{cmpA} provides the aggregated compound Latin cube:
\begin{equation}
\textbf{C}_\text{A} =
\begin{tabular}{llll}
$\fbox{$%
\begin{array}{rrrr}
0 & 1 & 2 & 3\\
1 & 0 & 3 & 2\\
2 & 3 & 0 & 1\\
3 & 2 & 1 & 0%
\end{array}%
$}$ & $\fbox{$%
\begin{array}{rrrr}
1 & 0 & 3 & 2\\
0 & 1 & 2 & 3\\
3 & 2 & 1 & 0\\
2 & 3 & 0 & 1%
\end{array}%
$}$ & $\fbox{$%
\begin{array}{rrrr}
2 & 3 & 0 & 1\\
3 & 2 & 1 & 0\\
0 & 1 & 2 & 3\\
1 & 0 & 3 & 2%
\end{array}%
$}$ & $\fbox{$%
\begin{array}{rrrr}
3 & 2 & 1 & 0\\
2 & 3 & 0 & 1\\
1 & 0 & 3 & 2\\
0 & 1 & 2 & 3%
\end{array}%
$}$ \\
\multicolumn{1}{c}{1$^{\text{st}}$ Layer} & \multicolumn{1}{c}{2$^{\text{nd}%
} $ Layer} & \multicolumn{1}{c}{3$^{\text{rd}}$ Layer} & \multicolumn{1}{c}{4$^{\text{th}}$ Layer}%
\end{tabular}
\label{latinCubeAgg}
\end{equation}
Using eq. \ref{cmpD} gives the dispersed Latin compound cube,
\begin{equation}
\textbf{C}_\text{D} =
\begin{tabular}{llll}
$\fbox{$%
\begin{array}{rrrr}
0 & 2 & 1 & 3\\
2 & 0 & 3 & 1\\
1 & 3 & 0 & 2\\
3 & 1 & 2 & 0%
\end{array}%
$}$ & $\fbox{$%
\begin{array}{rrrr}
2 & 0 & 3 & 1\\
0 & 2 & 1 & 3\\
3 & 1 & 2 & 0\\
1 & 3 & 0 & 2%
\end{array}%
$}$ & $\fbox{$%
\begin{array}{rrrr}
1 & 3 & 0 & 2\\
3 & 1 & 2 & 0\\
0 & 2 & 1 & 3\\
2 & 0 & 3 & 1%
\end{array}%
$}$ & $\fbox{$%
\begin{array}{rrrr}
3 & 1 & 2 & 0\\
1 & 3 & 0 & 2\\
2 & 0 & 3 & 1\\
0 & 2 & 1 & 3%
\end{array}%
$}$ \\
\multicolumn{1}{c}{1$^{\text{st}}$ Layer} & \multicolumn{1}{c}{2$^{\text{nd}%
} $ Layer} & \multicolumn{1}{c}{3$^{\text{rd}}$ Layer} & \multicolumn{1}{c}{4$^{\text{th}}$ Layer}%
\end{tabular}
\label{latinCubeDisp}
\end{equation}
For a cube, the rows and columns in each layer, as well as each column (sum of the elements through the layers at a given position in the cube) have identical line sums. The smallest order magic cube that can be written is $N=3$, of which there are four distinct cubes. The example shown above is straightforward to adapt to the magic case provided that $k=3$ for magic cubes since the elements in the cube range from $0..(N^3-1)$.

\end{appendices}

\end{document}